\documentclass{ifacconf}
\usepackage{natbib}        

\usepackage{savesym}
\savesymbol{AND} 
\usepackage{amsmath}
\usepackage{amssymb}
\usepackage{algorithm}
\usepackage{algorithmic}
\usepackage{graphicx}
\usepackage{nicefrac}

\newtheorem{proposition}{Proposition}[section]

\usepackage{dblfloatfix}

\newcommand{\cmt}[1]{}
\newcommand{\new}[1]{#1}

 \usepackage{xcolor}
 \usepackage[normalem]{ulem}
\newcommand{\remove}[1]{\cmt{{\color{red}{\sout{#1}}}}}
\newcommand{\repl}[2]{\remove{#1}\new{#2}}

\usepackage{xspace}

\newcommand\Tstrut{\rule{0pt}{2.3ex}}

\newcommand{\ie}{i.\,e.\xspace}
\newcommand{\eg}{e.\,g.\xspace}
\newcommand{\cf}{cf.\xspace}

\newcommand{\thisalg}{LNMS\xspace}
\newcommand{\moseq}{\ensuremath{\mathcal{M}}\xspace}
\newcommand{\feasr}{\ensuremath{\mathcal{X}^{\mathcal{M}}}\xspace} 
\begin{document}
\begin{frontmatter}

\title{Fast Non-Parametric Learning to Accelerate Mixed-Integer Programming for\\Hybrid Model Predictive Control\thanksref{footnoteinfo}}

\thanks[footnoteinfo]{Jia-Jie Zhu is supported by funding from the European Union’s Horizon 2020 research and innovation programme under the Marie Sk\l{}odowska-Curie grant agreement No 798321.}

\author[First]{Jia-Jie Zhu}
\author[Second]{Georg Martius}

\address[First]{Empirical Inference Department,\\Max Planck Institute for Intelligent Systems, T\"ubingen, Germany.\\(e-mail: jzhu@tuebingen.mpg.de)}
\address[Second]{Autonomous Learning Group,\\Max Planck Institute for Intelligent Systems, T\"ubingen, Germany.\\(e-mail: georg.martius@tuebingen.mpg.de)}

\begin{abstract}                
Today's fast linear algebra and numerical optimization tools have pushed the frontier of model predictive control (MPC) forward, to the efficient control of highly nonlinear and hybrid systems.
The field of hybrid MPC has demonstrated that exact optimal control law can be computed, \eg, by mixed-integer programming (MIP) under piecewise-affine (PWA) system models.
Despite the elegant theory, online solving hybrid MPC is still out of reach for many applications.
We aim to speed up MIP by combining geometric insights from hybrid MPC, a simple-yet-effective learning algorithm, and MIP warm start techniques.
Following a line of work in approximate explicit MPC, the proposed learning-control algorithm, LNMS, gains computational advantage over MIP at little cost and is straightforward for practitioners to implement.
\end{abstract}


\end{frontmatter}
\section{Introduction}


The presence of hybrid dynamical systems, defined as systems whose state evolution is governed by both continuous dynamics (flow) and discrete events (jump), is ubiquitous in physical systems.
Consider a legged robot that dynamically navigates the environment.
It must not only coordinate the motion of its joints but also decide when to make and break contact at its end-effectors.
This inherent coupling of \textit{discrete events} and \textit{continuous decision-making} challenges optimization-based control designs such as model predictive control (MPC).
One of the greatest difficulties is perhaps the combinatorial growth of computational complexity caused by mode switches.
\begin{figure}
  \centering
  \begin{tabular}{c@{\ }c}
    (a)&(b)\\
    \includegraphics[width=0.48\linewidth]{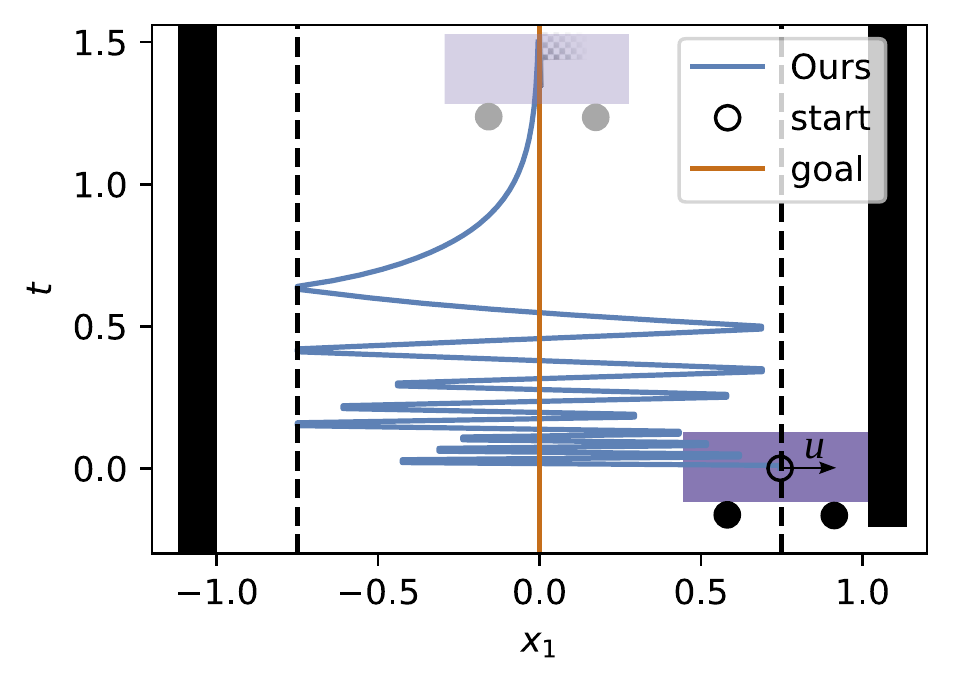}&
    \includegraphics[width=0.48\linewidth]{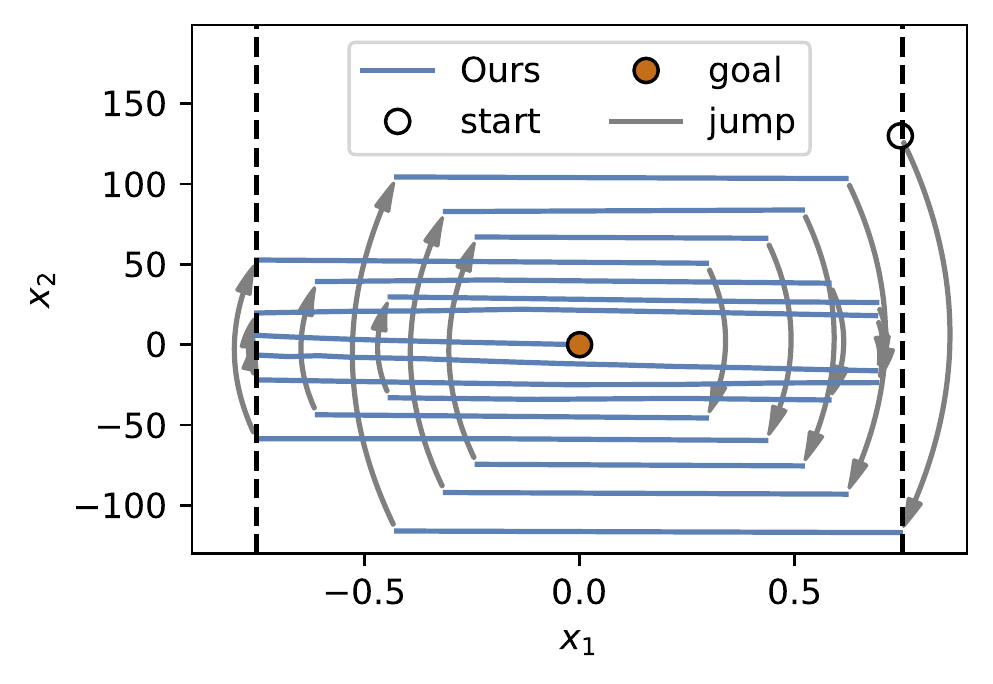}
  \end{tabular}
  \vspace{-1.2em}
  \caption{Illustrative example of Cart with two walls. (a) Cart starts at high velocity and bounces (discrete timesteps), note the vertical time axis. (b) same as in (a) but in state space ($x_1$:position and $x_2$: velocity). The gray lines with arrow denote jumps across switching modes. See Sec.~\ref{sec:exp} for details.}
  \label{fig:cart_running}
\end{figure}
This paper aims to answer an intuitive question: what information can be \textit{learned} when solving hybrid MPC problems over and over again?
Our \emph{main insight} is: The geometric structure of hybrid MPC solutions --- the mode sequences --- can be learned by a nonparametric learning algorithm from previously visited states.

Instead of storing exact solution structures, \eg polyhedral partitions, we store {samples of visited states} to greatly speed up hybrid MPC at little extra computational effort.


As an illustrative example, consider the dynamics of a cart in front of a wall, as shown in Figure \ref{fig:cart_running} (adapted from \cite{marcucci2019mixed}). It can be thought of as an actuated version of a bouncing ball, whose dynamics is continuous, except when it interacts with the wall.
%

During continuous dynamics periods, tools from continuous optimization such as sensitivity analysis,
that give rise to efficient \emph{online} methods \cite{ferreau2008online} or sequential quadratic programming (SQP)-based methods~\cite{diehl2002real},
can be used to control the system in a receding horizon fashion.
However, when contact happens, the structured solution of the sensitivity analysis loses validity.
%
%
In such cases, online hybrid MPC using mixed-integer programming (MIP) can be an effective way of making optimal discrete decisions (\eg whether or not to willingly impact the wall). However, solving MIP online is computationally heavy.


We hope to contribute in the following aspects:
\begin{itemize}
	\item We propose a simple-yet-effective fast learning algorithm to warm-start mixed-integer programs for MPC in PWA hybrid systems. As learning progresses, it greatly reduces the computational cost.
	\item The properties of the proposed algorithm are exploited to allow i)~post-processing of solutions improving towards optimality and early-termination of MIP; ii)~straightforward implementation for practitioners.
\end{itemize}


\textbf{Notation}
In this paper, $x\in R^d$ denotes a column vector of dimension $d$ and $x^\top$ its transpose. $N$ denotes the horizon of an optimal control problem (OCP) and $n$ the number of samples of the learning algorithm. $x_p$ is a parameter (\eg current state estimation) of an optimization problem.
$e_i$ denotes the one-hot vector of suitable dimensions where the $i$-th element is $1$ and rest are zeros.
A sequence of modes of a discrete dynamical system is denoted as $\moseq=\{m_t\}^N_{t=1}$ where each $m_t = (\mu_t^1, \mu_t^2,\dots,\mu_t^{n_M})$ is a system mode at time $t$ and $\mu_t^i$ a binary variable.
A feasible solution of an optimization problem is one at which all constraints are satisfied.

\section{Preliminaries}
\label{sec:prelim}
\subsection{Model predictive control for piecewise affine hybrid systems}
In this paper, we consider a discrete-time optimal control problem (OCP)\remove{in the following optimization problem}.
The goal is to compute a sequence of control actions $\{u_t\}_{t=1}^N$ to steer the system state $x$ to the origin.
\begin{equation}
\label{eq:ocp}
\begin{aligned}
& \underset{u_t, t=0,\ldots,N-1}{ \text{minimize}} & & \sum_{t=0}^{N-1} \left(x_t^\top Q x_t + u_t^\top  R u_t\right) + x_N^{\top}P x^{ }_N,\\
& \quad\text{subject to} & & x_{t+1} = f(x_t, u_t),\ t=0,\ldots,N-1,\\
& & & h_t(x_t, u_t) \leq 0,\ t=0,\ldots,N-1,\\
&&& x_0 = x_p
,
\end{aligned}
\end{equation}
where $x_p$ is the current state, $x_t^\top Q x_t + u_t^\top  R u_t$ is the stage cost and $x_N^\top  P x^{ }_N$ the terminal cost. $f$ is the dynamics that describes the evolution of the system over time and $h_t$'s are the constraints (\eg, bounds on control input $u_t$).

Model predictive control (MPC)~\cite{richalet1978model, rawlings2009model} solves OCP~\eqref{eq:ocp} at every sampling time and implement the first control input $u_0$.

Consider the running example of the simple cart and wall in Fig.~\ref{fig:cart_running} again. In this case, there is no single dynamics function $f$ that governs on the whole state-space --- this is a hybrid system. One approach to solve hybrid MPC is to formulate the dynamics constraints using piecewise affine (PWA) formulations (\cite{sontag1981nonlinear}).
PWA models can describe general nonlinear dynamics while enjoying the advantage of having a wide range of linear control tools,
\begin{equation}
\label{eq:pwa}
\begin{aligned}
& x_{t+1} = A^{i}x_t + B^{i}u_t, \text{ for }x\in \mathcal{C}_i, i=1\dots,n_M.
\end{aligned}
\end{equation}
Mathematically, solving OCP~\eqref{eq:ocp} under PWA dynamics~\eqref{eq:pwa} is typically formulated as a mixed-integer program (MIP), due to the presence of both continuous and discrete variables.
We consider the big-M formulation of MIP:
\begin{equation}
\label{eq:bigm}
\begin{aligned}
|x_{t+1} - A^{i}x_t - B^{i}u_t| &\leq (1-\mu^{i}_t)  M\\
h(x_t, u_t)&\leq (1-\mu^{i}_t)  M\\
\sum_i \mu^{i}_t &= 1, \ \mu^{i}_t \in \{0,1\},\quad \forall i,t,
\end{aligned}
\end{equation}
where $\mu^{i}_t, i=1,2,\dots n_M$ are the auxiliary integer variables.
At time $t$, it is easy to see that the system is in mode $i$: $m_t=e_i$ if and only if $ \mu^{i}_t=1$. A mode sequence $\moseq=\{m_t\}^N_{t=1}$ thus corresponds to a set of integer decision variables of the hybrid OCP \eqref{eq:ocp},\eqref{eq:bigm}.

We use \feasr to denote the feasible region (not to be confused with the critical region)
of an OCP given\remove{known} \moseq: the set of parameters $x_p$ where a mode sequence \moseq  is feasible (\ie, OCP \eqref{eq:ocp}\eqref{eq:bigm} has a solution with this fixed \moseq).
It is easy to see that one mode sequence may not be feasible for the whole state space in hybrid systems.
For example, consider a PWA control law of the cart example in Fig.~\ref{fig:voronoiknn2000} (a), the mode sequence associated with the region (colored red) in the upper right corner is not feasible once we cross the boundary of this region.
This corresponds to a different set of integer solutions to OCP \eqref{eq:ocp}\eqref{eq:bigm}.
The geometric property of feasible regions can be summarized as follows:
i)~Feasible regions are convex polyhedra. ii)~Each feasible region corresponds to an integer solution to the hybrid OCP\eqref{eq:ocp}\eqref{eq:bigm} (and hence a mode sequence). iii)~They may or may not overlap.

One insight commonly exploited is that
\repl{if we fix integer variables in \eqref{eq:bigm} to a \emph{feasible mode sequence} \moseq, the hybrid OCP \eqref{eq:ocp}\eqref{eq:bigm} becomes a continuous (convex) optimization problem and can be solved by extremely efficient algorithms.}{
the hybrid OCP \eqref{eq:ocp}\eqref{eq:bigm} becomes a continuous (convex) optimization problem
 if the integer variables in \eqref{eq:bigm} are fixed  to a \emph{feasible mode sequence} \moseq.
Thus, the problem can be solved by extremely efficient algorithms.
} Hence, if we store the feasible region \feasr and associated \moseq offline, during the online runtime of MPC, we only need to \emph{look up which \feasr the current state $x_p$ belongs to and fix the mode sequence to associated \moseq}.


\subsection{Voronoi tessellation and nearest neighbor classification}
One of the simplest yet effective machine learning algorithms is nearest neighbor (NN) classification [\cite{cover1967nearest}].
It stores all data in terms of input-output pairs. For a new data point, the
 nearest neighbor is retrieved and its associated output is returned.
 Let $D$ be the set of data points with entries $x_i \in \mathbb R^d$.
 The runtime complexity depends on the underlying storage structure: with tree structures [\cite{Omohundro89fiveballtree}]
  the retrieval is $O(d\log(|D|))$ (while building the indexing structure requires $O(d |D| \log(|D|))$ operations.).
Interestingly, the explicit MPC approach in \cite{jones2006logarithmic} uses Voronoi diagram as a way to reduce the complexity of online execution.
\repl{They report online look-up complexity logarithmic in the number of regions. In comparison, we use Voronoi tessellation but in a sampling setting which scales logarithmically with the  number of samples.}{
Their online look-up complexity is logarithmic in number of regions
whereas our setting scales logarithmically with the number of samples.
This gives us a flexible trade-off: we can choose to store fewer samples with faster look-up time but more likely to need to solve MIP, or to store more samples for the opposite, or somewhere in between.
}

\begin{figure}
	\centering
	\begin{tabular}{c@{}c@{}c@{}c}
		(a) PWA control law & (b) ground truth\\
		\includegraphics[width=0.45\linewidth]{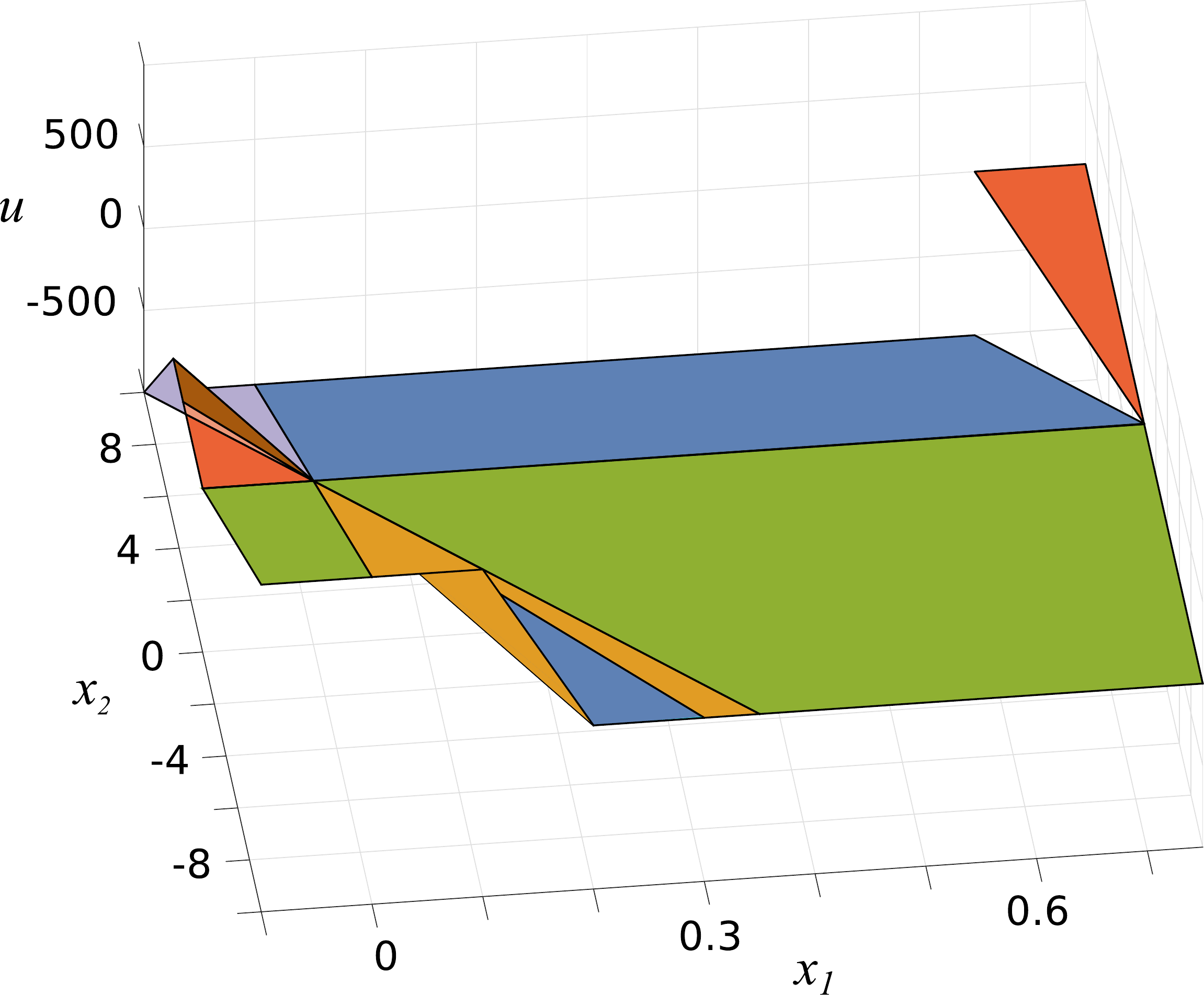}
                 &\includegraphics[width=0.45\linewidth]{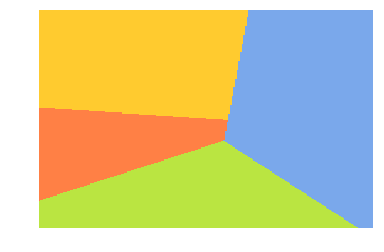} \\
              (c) 100 samples&(d) 2000 samples\\
		\includegraphics[width=0.45\linewidth]{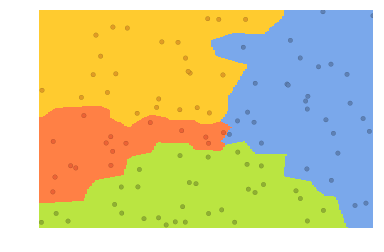}
		&\includegraphics[width=0.45\linewidth]{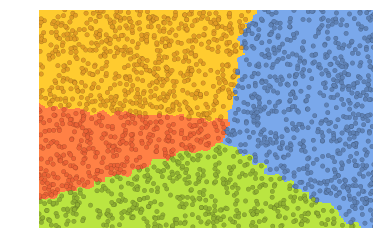}
	\end{tabular}
  \vspace*{-0.4em}
	\caption{(a) piecewise affine control law in hybrid systems (cart example). The z-axis is the control input while x,y-axes are states. Notice the control law is discontinuous.
\remove{Note such plots are typically visualized in constrained linear control studies, e.g. \mbox{\cite{mayne2003optimal}}.}
    (b) Ground truth NN partition of regions.
    (c,d) Evolution of Voronoi regions induced by NN with increasing numbers of sampling points.
    }
    	\label{fig:voronoiknn2000}
\end{figure}

Fig.~\ref{fig:voronoiknn2000} (b,c,d) is an intuitive example of NN classification resulting in the Voronoi diagram approximating the true region partition. Such approximation is the \emph{key motivation} to our method:
we directly store sampled points rather than keeping track of polyhedral regions in Fig.~\ref{fig:voronoiknn2000} (a).
Sampled points offer a straightforward alternative to solving multi-parametric (quadratic) programs  or performing polyhedron operations.

\emph{Remark}
One key feature of NN is that it belongs to the class of ``\textit{lazy learning}'' algorithms; it does not require an additional training process such as neural networks do. This important aspect facilitates the \textit{fast-learning} capability of the learning-controller we shall propose next.

\section{Method}
\label{sec:method}
\subsection{Approximate hybrid MPC by learning from nearest neighbor mode sequences}
\label{sec:knn_mpc}
Drawing form geometric insights  discussed in Sec.~\ref{sec:prelim}, we propose a concise approximation approach based on the nearest neighbor (NN) rule.
Our \emph{main idea} is simple: use Voronoi tessellation induced by the
NN classifier to approximate feasible regions for the hybrid MPC control law.
The data for the NN classifier is given by the mode sequence \moseq based on the \emph{sampled} states.

The learning-control algorithm works as follows. We query the NN classier for the mode sequence \moseq for a given state $x_p$ based on which feasible region it falls in (a simple NN look-up). We then use \moseq's corresponding integer solution to warm-start the MI(Q)P \eqref{eq:ocp},\eqref{eq:bigm}. If the modes are feasible, the computation is reduced to a quadratic program that can be solved with state-of-the-art solvers in the order of microseconds, \eg \cite{houska2011acado}.
The algorithm is explained in Alg.~\ref{alg:hmpc}. For the rest of the paper, we refer to the proposed method as \emph{approximate MPC by learning from nearest-neighbor mode sequences} (\thisalg).

\begin{algorithm}[tbp]
	\caption{\thisalg: Online Approximate MPC with Nearest Neighbor Learning}
	\label{alg:hmpc}
	\begin{algorithmic}[1]
		\STATE Given: sample initial state $x(0)\sim P_0$ where $P_0$ is the initial state distribution. Optional: an initial dataset $\mathcal{D} = \mathcal{D}_0$ of sampled points.
		\LOOP
		\STATE Get current state estimation $x_p$ .
		\STATE Query the nearest neighbor classifier (with dataset $\mathcal{D}$) for the mode sequence \moseq$=\{m_i\}^N_{i=1} $
		\STATE Solve hybrid OCP~\eqref{eq:ocp}\eqref{eq:bigm} with integer variables $\{m_i\}^N_{i=1}$ as warm-start solution. Terminate once a feasible solution is found.
		\STATE Add the $(x_p, \mathcal{M^*})$-pair to the dataset $\mathcal{D}$, where $\mathcal{M^*}$ is the integer solution obtained last step.
		\STATE Apply the first control solution $u_0^*$ to the system.
		\ENDLOOP
	\end{algorithmic}
\end{algorithm}





\subsection{Theoretical properties}
\label{sec:theory}
We present the theoretical analysis for Alg.~\ref{alg:hmpc}.
Typical MPC analysis concerns feasibility (whether the OCP has a solution) and stability (whether the closed-loop system can be bounded around a set-point or within a set\remove{ in the future}).
For (exact) hybrid MPC, many properties are inherited from nominal MPC, \ie \emph{recursive feasibility (and therefore asymptotic stability) is guaranteed} by choosing the appropriate terminal cost (Lyapunov function) and terminal constraint (control invariant set).
We refer interested readers to \cite{borrelli2017predictive}, Section~17.8.1 for detailed theoretical analysis regarding HMPC stability.
As Alg.~\ref{alg:hmpc} executes the regular hybrid MPC as a warm-start improvement in the case of predicted modes being infeasible, it inherits the feasibility (hence stability) properties from regular hybrid MPC. We will focus on the computational aspect, i.e., the speed-up over exact hybrid MPC.

\begin{proposition}
  Given dataset $\mathcal{D}$ of size $n$ obtained by Alg.~\ref{alg:hmpc}, let $P_n^{\text{MIP}}$ denote the probability of hybrid MPC, whose initial state $x(0)\sim P_0$, having to execute MIP solver in Step~5 of Alg.~\ref{alg:hmpc}. Then,
$$P_n^{\text{MIP}}\to 0 \text{  as  } n\to \infty.$$
\end{proposition}
\emph{Proof}
For every feasible state, we consider the MIQP solver to be an oracle. It achieves the minimum possible
classification error rate $P^* = 0$ as we assume the solver is always able to return feasible mode sequences for feasible states.
Let $P_n(\text{error})$ denote the probability that the NN classifier in Algorithm~\ref{alg:hmpc} predicts an infeasible mode sequence (classification error).
Use the NN-convergence bound $P^* \leq \lim_{n\to \infty}P_n(\text{error}) \leq P^* (2 - \frac{c}{c-1}P^*)=0$, where $c$ is the number of classes (c.f.\ \cite{cover1967nearest}). Since we only execute MIP solver if the mode sequence is infeasible, $P_n^{\text{MIP}}=P_n(\text{error})$. The
conclusion follows.$\hfill\square$

This proposition implies that, with an increasing amount of
samples, we get a faster and faster controller.
This is empirically validated in the experiment section.
%
%

\subsection{Improving optimality of sampled mode sequences by warm-starting MIP}
\label{sec:improve}

As an approximate MPC algorithm, \thisalg only aims to produce feasible instead of optimal mode sequence prediction.
In practice, one may also choose to early terminate the MIP at a feasible solution since the most costly computation is to produce a tight dual bound to certify optimality --- a good solution may be available much sooner [cf. \cite{bertsimas2016best}]. Both those two sources contribute to suboptimal mode sequences for sampled points.

However, using the instance-based nature of the NN classifier, we can improve on the resulting controller by simply ``relabeling'' the stored samples via warm-starting techniques of MIP.

Intuitively, this process picks a point from the stored samples and feed its mode sequence as warm-start to the MIP solver, see Alg.~\ref{alg:imprv}). As it is already feasible, the solver will always return a mode sequence that \repl{is more optimal}{has the same or lower cost}.
Hence we have the intuitive results in Proposition~\ref{thm:imprv} (omitting the straightforward proof). We provide examples of this process in Sec.~\ref{sec:exp} (Fig.~\ref{fig:e1_static}(d), Fig.~\ref{fig:twowall}(c)).

\begin{algorithm}
	\caption{\thisalg: Solution optimality improvement}
	\label{alg:imprv}
	\begin{algorithmic}[1]
		\LOOP
		\STATE Pick a state-label $(x_0, \mathcal{M})$-pair from dataset $\mathcal{D}$

		\STATE Solve OCP~\eqref{eq:ocp} with integer variables \moseq$=\{m_i\}^N_{i=1}$ as a warm-start incumbent solution. Optionally, early terminate after the computational budget reached.
		\STATE Relabel this sample and add the new $(x_0, \mathcal{M*})$-pair to $\mathcal{D}$, where $\mathcal{M^*}$ is the integer solution from step 3..,
		\ENDLOOP
	\end{algorithmic}
\end{algorithm}

\begin{proposition}
	\label{thm:imprv}
	If the re-labeling using Algorithm~\ref{alg:imprv} is done to full MIP optimality for all stored data points, then Algorithm~\ref{alg:hmpc}, with dataset $\mathcal{D}$ and initial state $x(0)\sim P_0$,  recovers the exact optimal solution of hybrid MPC as the number of samples $|\mathcal{D}|\to \infty$.
\end{proposition}

\section{Numerical Experiments}
\label{sec:exp}

\subsubsection{Experiment setup}
We implement the MPC controller in Python with Gurobi as the optimization solver.
In NN learning, we use a simple weighted Euclidean distance. \repl{ (\eg in 2D vector space, $\|p\|_w = \sqrt{w_1 p_1^2 + w_2 p_2^2}$. If $w_1=w_2=1$, this reduces to Euclidean distance). }
We consider two examples using a cart with one and two walls, see Fig.~\ref{fig:cart_running}(a) (in the one wall case the left wall is missing), and a pendulum with an elastic wall, see~Fig.~\ref{fig:e3}(a).
\subsection{Example 1. Cart-wall contact dynamics}

The dynamics equation for the cart-wall example in Fig.~\ref{fig:cart_running} is described by the following PWA formulation:
\begin{equation}
\label{eq:db_int}
\begin{aligned}
\begin{cases}
x_1^+ = x_1 +  x_2\ \Delta t, \quad x_2^+ = x_2 + \dfrac{u}{m}\ \Delta t,  & \text{if } x\in \mathcal{C}_1\\
x_1^+ = x_1, \quad x_2^+ = -\epsilon x_2 , & \text{if } x\in \mathcal{C}_2
\end{cases}
\end{aligned}
\end{equation}
where $m$ is the cart mass (set to 1.0) and $\epsilon$ is the coefficient of restitution (set to 0.9).
It could be thought of an actuated version of a bouncing ball --- a classic hybrid system.
 $\mathcal{C}_1 = \{x_1 +  x_2\ \Delta t < x_{\text{wall}}\}$ denotes the state space where the dynamics is the double integrator and
 $\mathcal{C}_2 = \{x_1 +  x_2\ \Delta t \ge x_{\text{wall}}\}$ denotes the state space where the contact with the wall happens. In our case $x_{\text{wall}} = 0.75$ and the discretization is set to $\Delta t = 0.01$ for all our experiments.
The PWA dynamics for the cart with two walls is a straightforward extension.

\begin{figure}
  \centering\scriptsize
  \begin{tabular}{cc}
    (a) 15 samples (2 region)&(b) 100 samples (3 regions)\\
     \includegraphics[height=0.2\textwidth]{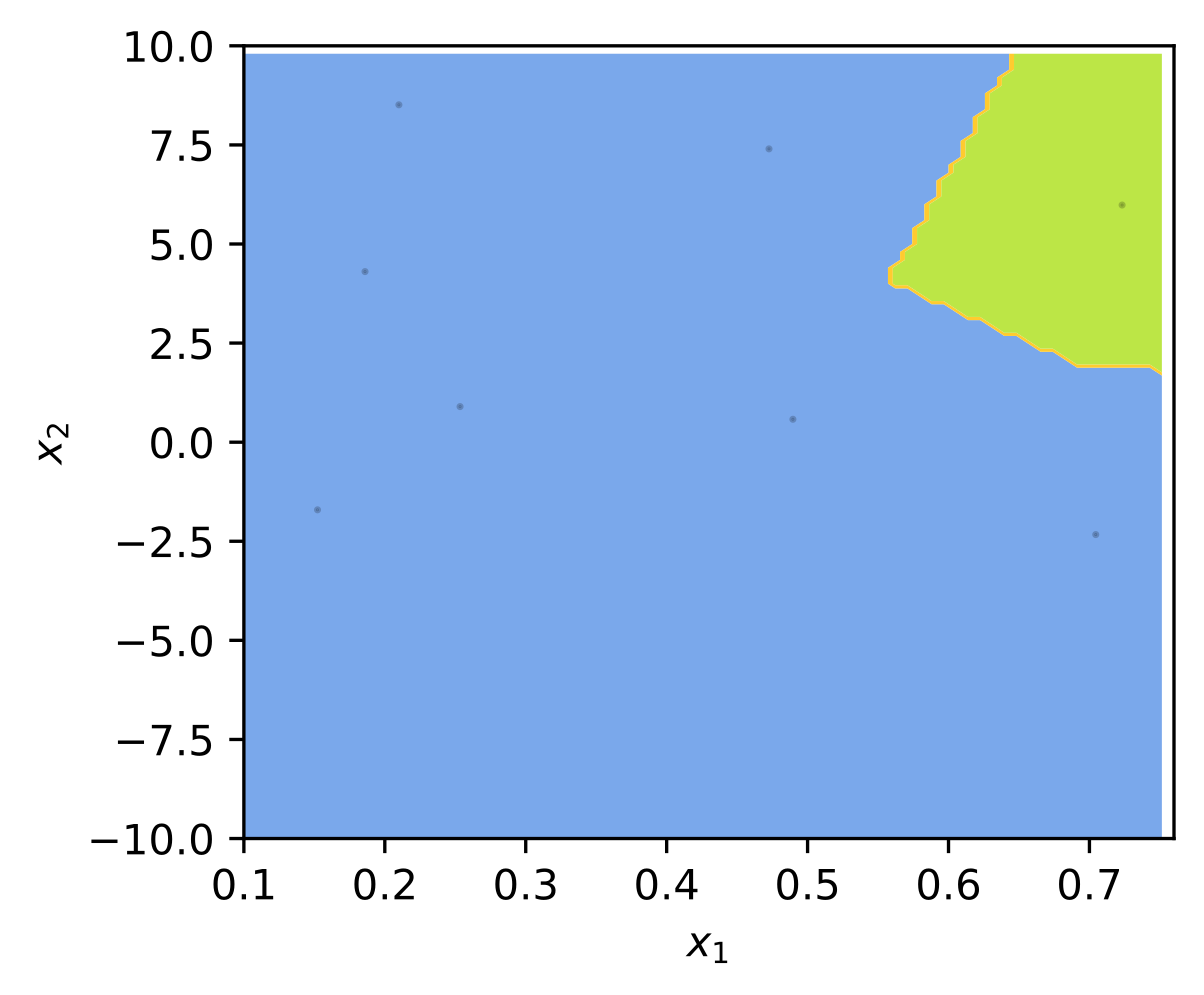}&
    \includegraphics[height=0.2\textwidth]{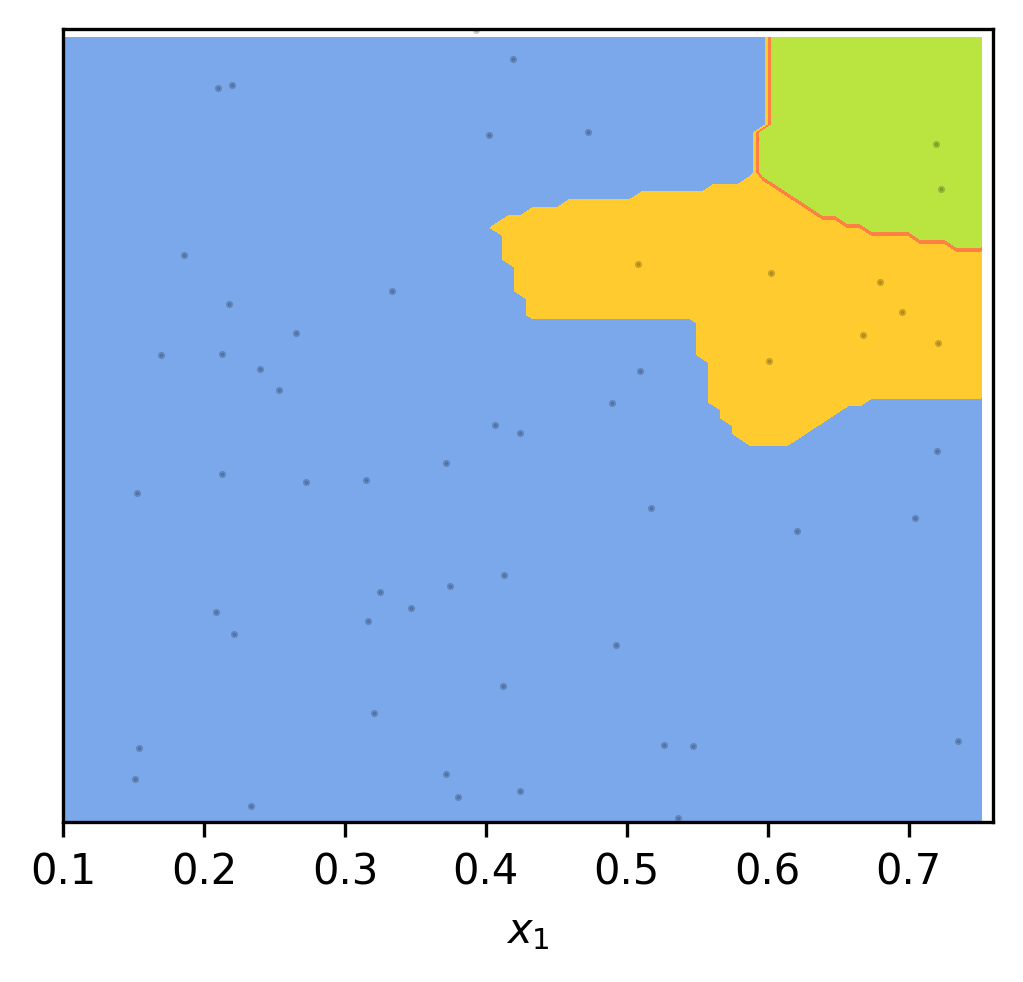}\\
    (c) 1000 samples (3 reg.)&(d) improved (4 regions)\\
   \includegraphics[height=0.2\textwidth]{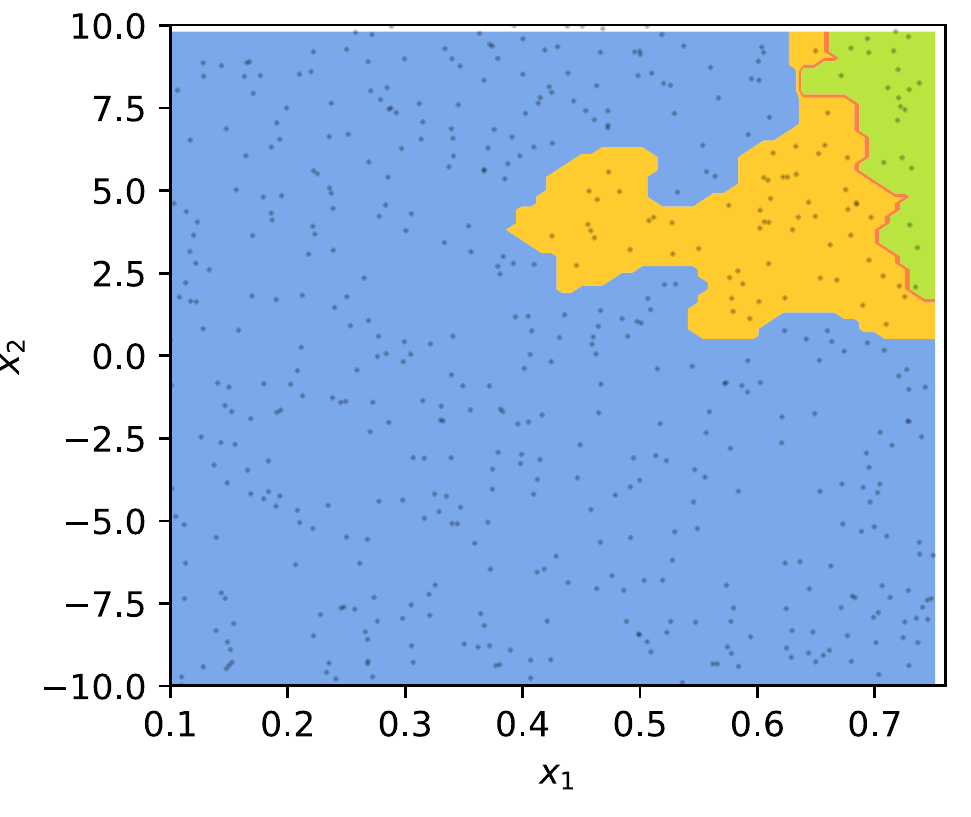}&
   \includegraphics[height=0.2\textwidth]{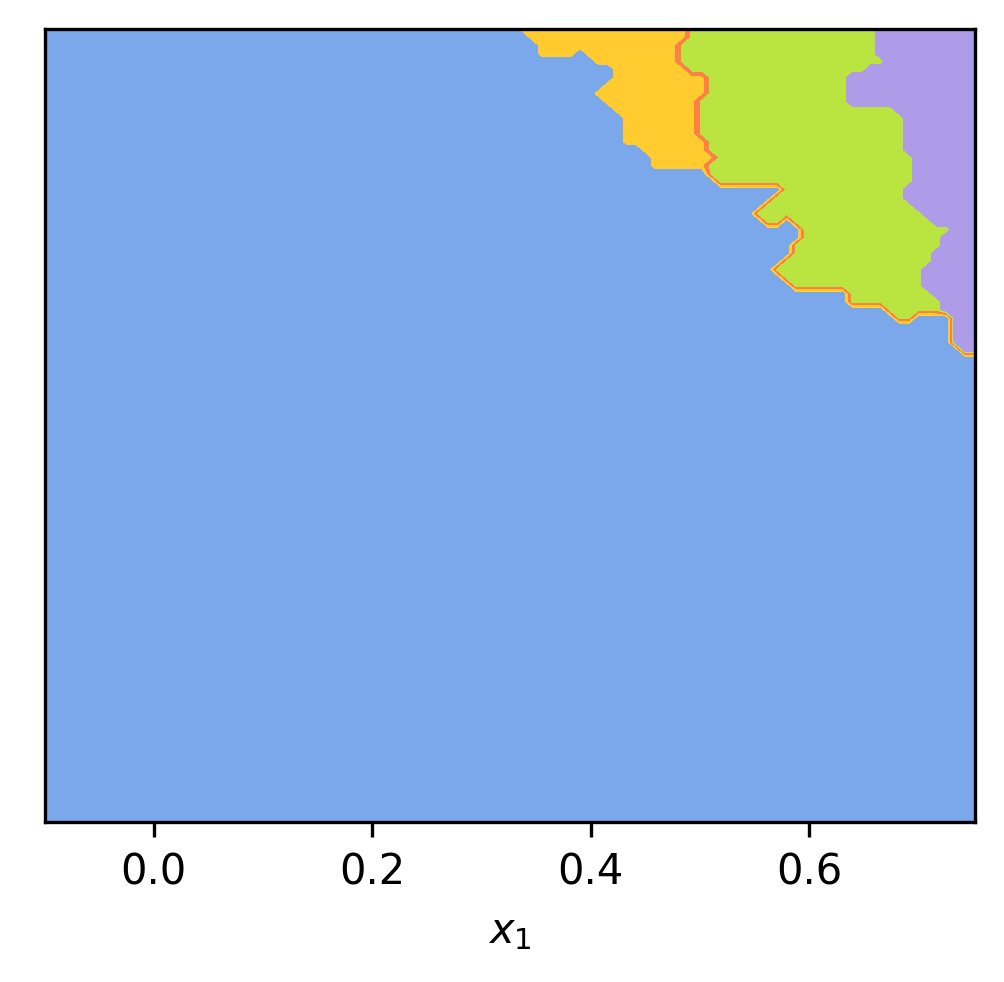}
  \end{tabular}
  \vspace*{-1em}
  \caption{PWA control law resulting from \thisalg in Alg.~\ref{alg:hmpc} associated with different sample sizes shown in (a-c).
    Panel (d) displays the optimal control law after applying the improvement scheme in Sec.~\ref{sec:improve}.
  }
  \label{fig:e1_static}
\end{figure}

We synthesize MPC with horizon $N=10$ (no early termination of the MIP).
The cost weights in Eqn.~\eqref{eq:ocp} are
$Q=I_2, R=0.001$. $S=\beta \cdot K$ where $K$ is the solution to the algebraic Riccati equation for the system in mode $\mathcal{C}_1$, we choose $\beta=1000$ for faster convergence behavior to the attractor but this is not a crucial choice.
In the first experiment, we use terminal costs instead of terminal sets following common practice in MPC applications (\cf \cite{rawlings2009model}).



In this experiment, we start from initial states $x(0)$ randomly sampled within the region $[0.1, 0.75]\times [-10,10]$. For each initial state, the hybrid OCP~\eqref{eq:ocp}\eqref{eq:bigm} is solved by \thisalg in Alg.~\ref{alg:hmpc}.
As a result of the non-parametric learning Alg.~\ref{alg:hmpc}, we obtained a set of samples $\mathcal D$. Those samples store the solution structure (feasible mode sequences) of \thisalg and are the key to speeding up MIP. They are displayed in Fig.~\ref{fig:e1_static} (a)-(c) (gray dots) with different sample size.

\begin{figure*}[b]
  \centering\footnotesize
  \begin{tabular}{ccc}
    (a) 5000 samples&(b) zoomed center of (a) &(c) improved version of (b)\\
    493 regions&153 regions &33 regions\\
    \includegraphics[width=0.3\textwidth]{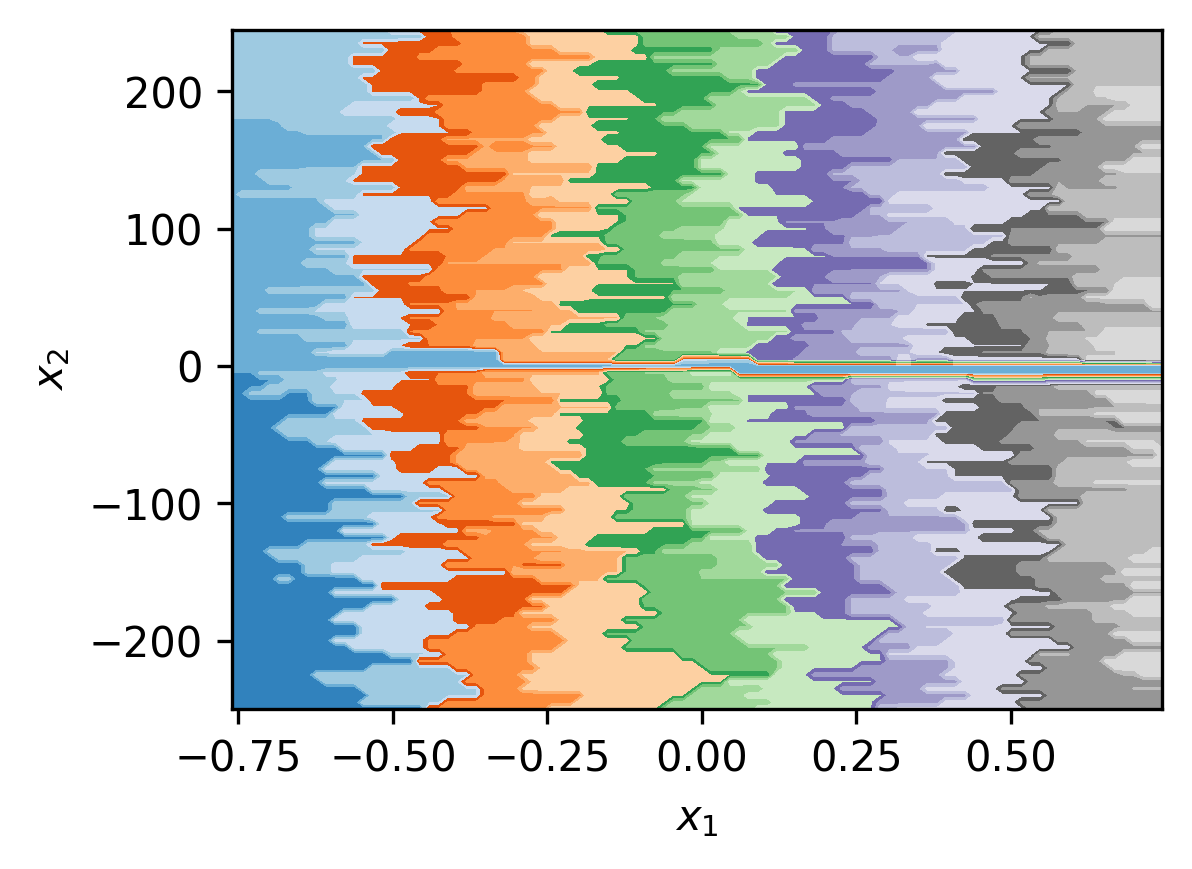}&
    \includegraphics[width=0.3\textwidth]{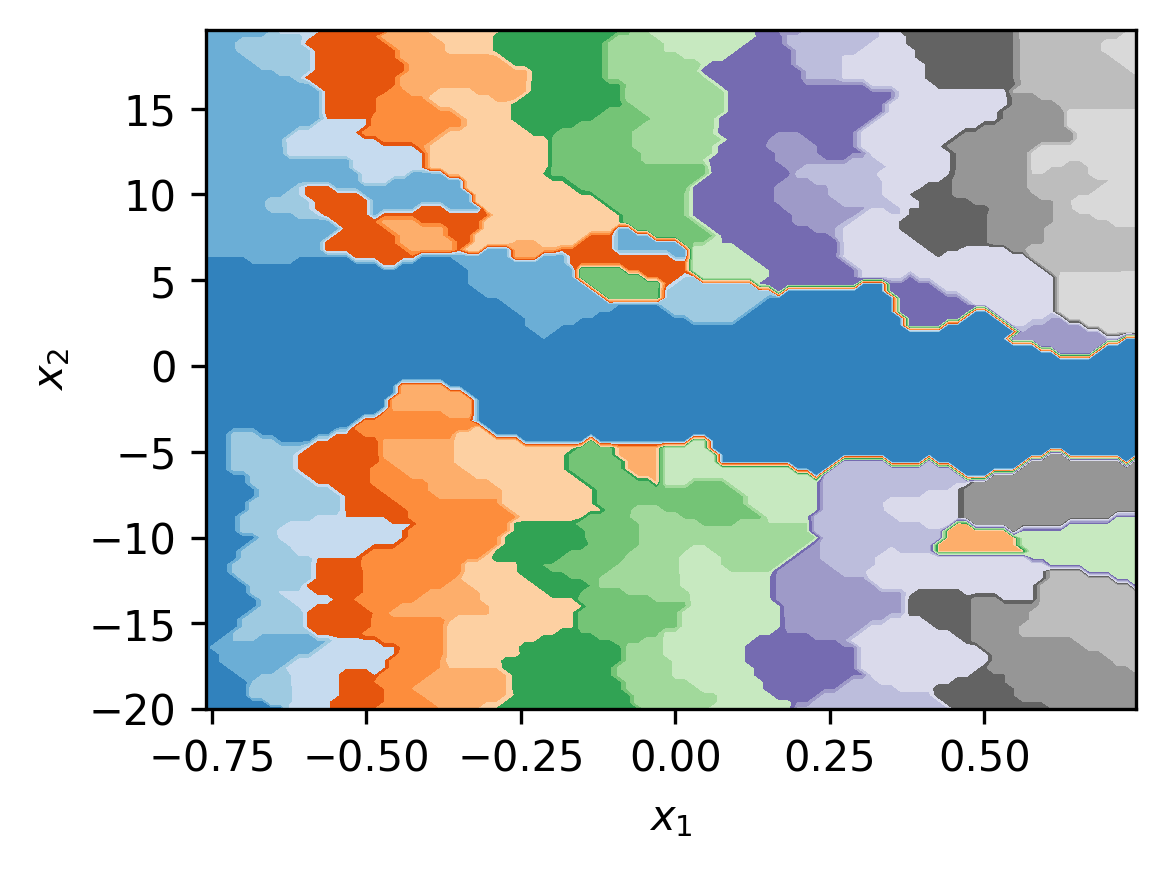}&
    \includegraphics[width=0.3\textwidth]{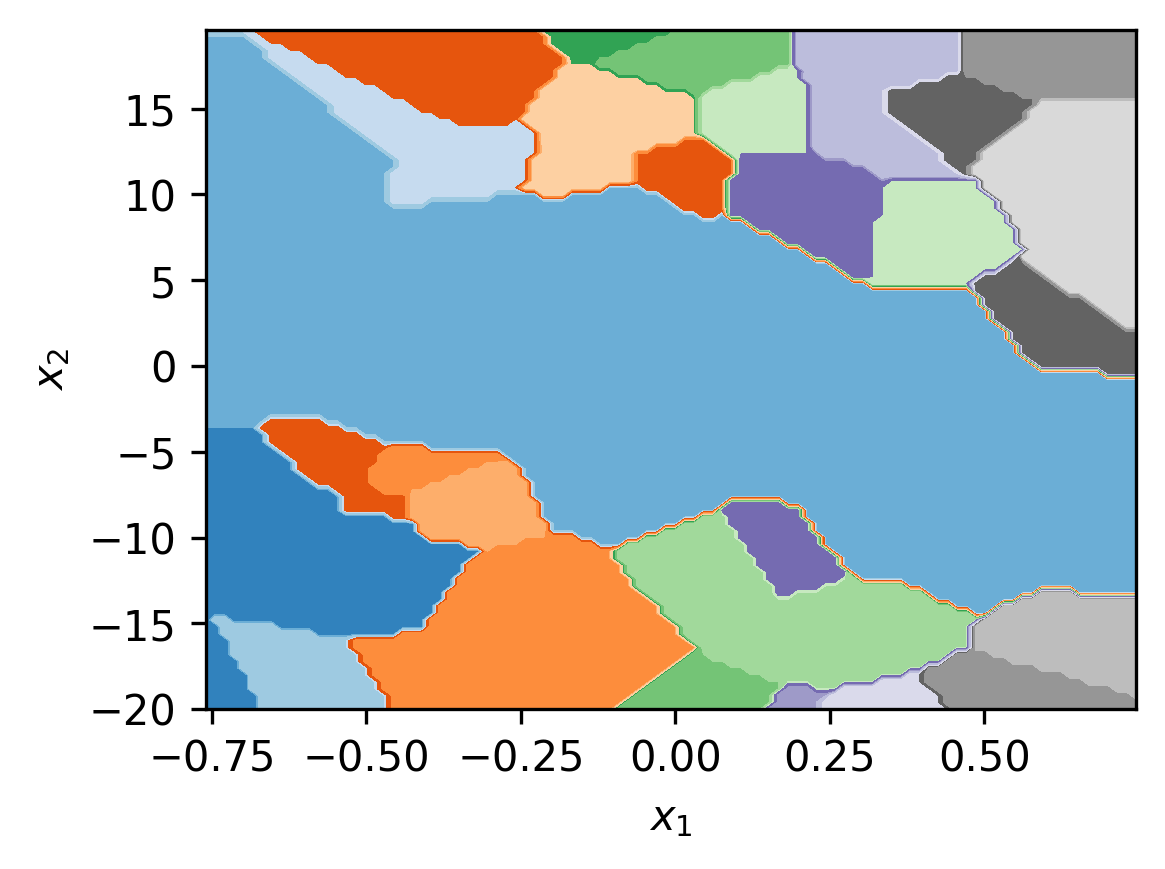}
  \end{tabular}
  \vspace*{-1em}
  \caption{Control law of \thisalg in the cart with two walls.
    (a) PWA control law of Alg.~\ref{alg:hmpc} with early termination of 5.0s. Different color patches denote different affine control policies.
    (b) zoomed view on the center area of (a).
    (c) after applying the warm-start relabeling method using 5s improvement proposed in Sec.~\ref{sec:improve}.
  }
  \label{fig:twowall}
\end{figure*}

\begin{figure}
  \centering
  \includegraphics[width=0.6 \linewidth]{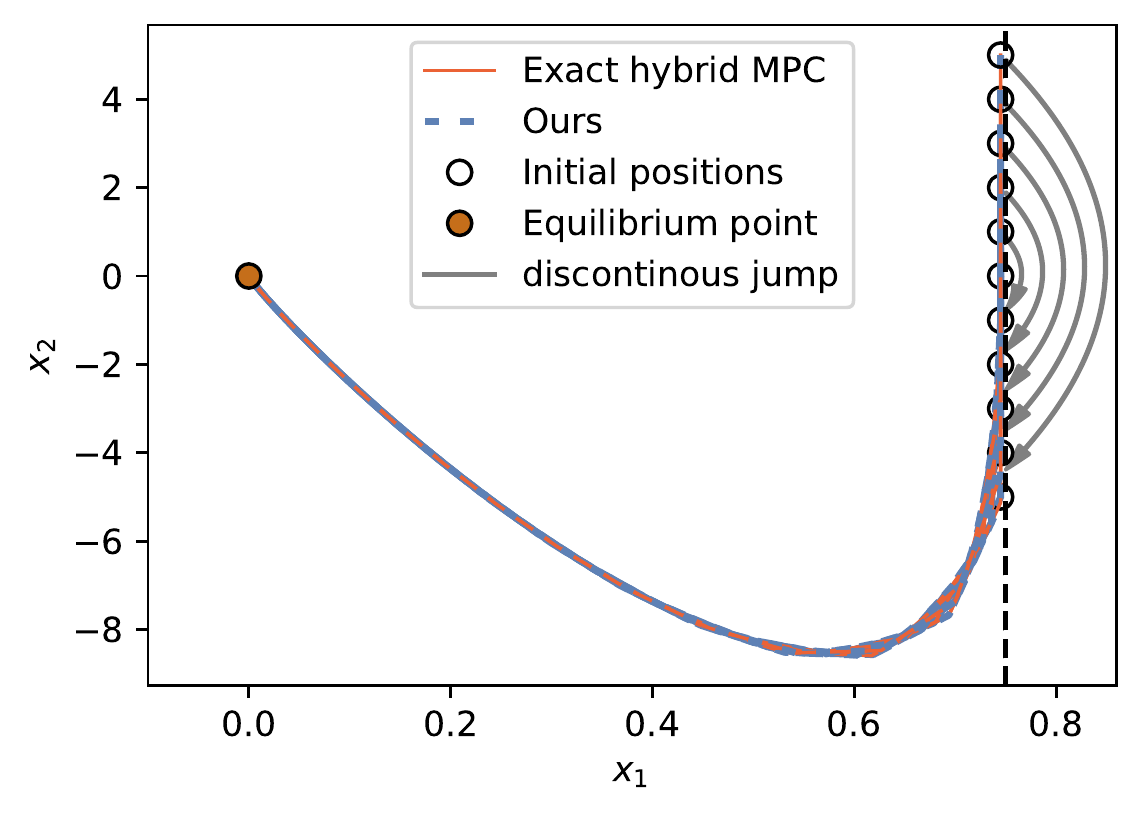}
  \vspace*{-1.5em}
  \caption{Example trajectory for the cart-1-wall.}
  \label{fig:ex1_traj}
\end{figure}

We visualize the resulting PWA control law in Fig.~\ref{fig:e1_static} (a)-(c). Each color patch indicates a region with different affine control laws.
Notably, Fig.~\ref{fig:e1_static} (a)-(c) shows the evolution of the PWA control law as the learning-controller gathers more samples. We refer to our earlier discussion on Fig.~\ref{fig:voronoiknn2000} for how to interpret the visualization.
We plot the resulting closed-loop state trajectories in Fig.~\ref{fig:ex1_traj} (with different initial states).
The trajectories obtained by \thisalg are identical to the exact hybrid MPC control law.


\subsubsection{Scaling up to complex control law structure}
It is known that constrained control law with numerous mode switches has a complex structure. To demonstrate that \thisalg scales to such cases, we consider the cart example with two walls, higher initial velocities, and an MPC prediction horizon of $N=25$.
Furthermore, we impose the input constraint ($|u| \leq 10$).
The resulting control law is significantly more complex than the previous example due to a large number of mode sequences and constraints.

A direct consequence is that the MIP is challenging to solve to full optimality.
Therefore, early termination is necessary --- we use a 5\,sec early-termination threshold.
The HMPC control law is shown in Fig.~\ref{fig:twowall}(a,b).
We observe a more complex structure of the control law --- 493 regions of different \repl{PWA control law as a result}{mode sequences correspond to different PWA control laws}.

\subsubsection{Improving solution optimality with Alg.~\ref{alg:imprv}}
As discussed, Alg.~\ref{alg:hmpc} seeks a feasible solution to speed up and warm-start MIP. The resulting solution might be a sub-optimal control law.
However, due to the non-parametric nature of our learning algorithm, we can simply post-process the stored sample points to improve the optimality of the learning-controller.
We apply the optimality-improvement scheme in Alg.~\ref{alg:imprv} to both the one-wall and two-wall cases.
The resulting \textit{optimal} control law is shown in Fig.~\ref{fig:e1_static}(d) for one-wall and Fig.~\ref{fig:twowall}(c) for two-walls.
Interestingly, we observe that the control law in Fig.~\ref{fig:twowall}(c) is simplified due to merging via comparing optimal cost-to-go of different PWA control laws. This is similar to the procedures in merging piecewise partitions in explicit MPC (cf. \cite{borrelli2017predictive}).

\begin{figure}
  \centering
  \begin{tabular}{c@{\ }c}
  (a) environment   &(b) partition (1000 samples)\\
  \includegraphics[width=0.3\linewidth, trim=0 -40 0 0]{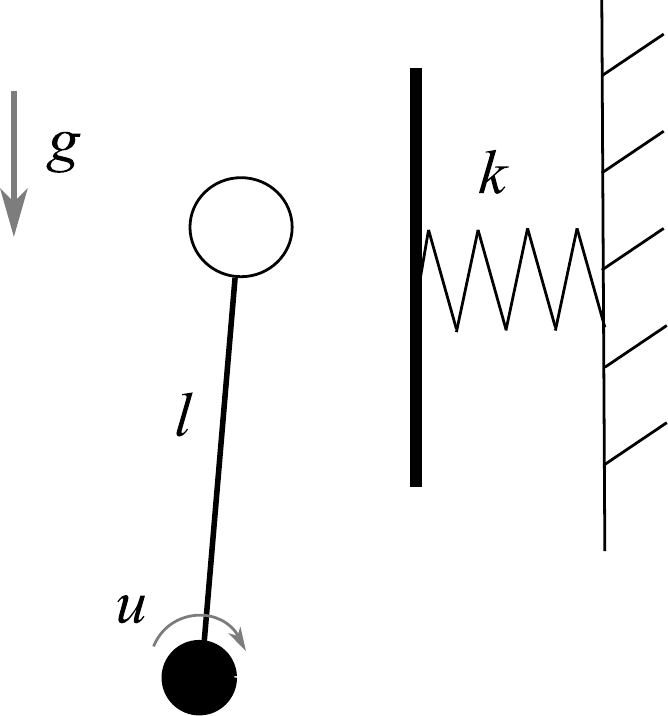}&
  \includegraphics[width=0.66\linewidth]{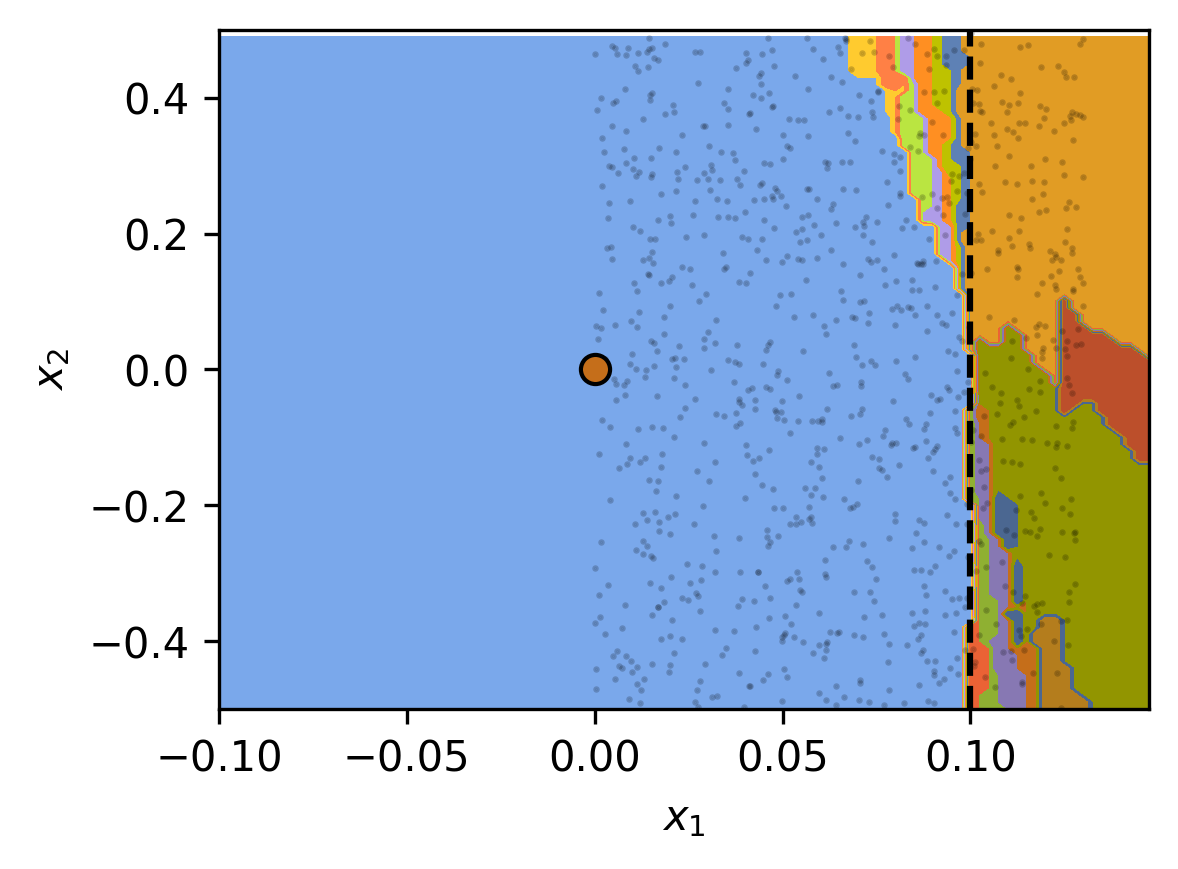}
  \end{tabular}
  \vspace{-.5em}
\caption{
  Pendulum with an elastic wall. (a) Illustration of the environment.
  (b) Partitioning from Alg.~\ref{alg:hmpc} resulting in 18 regions.
}
\label{fig:e3}
\end{figure}


\subsection{Example 2. Elastic wall with variable duration contact}
\label{sec:elastic_wall}
The elastic pendulum environment, see Fig.~\ref{fig:e3}(a), is adapted from~\cite{marcucci2017approximate}.
The dynamics is given by
\begin{equation}
\dot x =
\begin{cases}
A_1 x + B_1 u, \quad \text{if} \quad (x, u) \in \mathcal D_1, \\
A_2 x + B_2 u + c_2, \quad \text{if} \quad (x, u) \in \mathcal D_2,
\end{cases}
\end{equation}
with
\begin{align*}
A_1 &=
\begin{bmatrix}
0 & 1 \\ \nicefrac{g}{l} & 0
\end{bmatrix},
\quad
B_1 =
\begin{bmatrix}
0 \\ \nicefrac{1}{ml^2}
\end{bmatrix},\\
A_2 &=
\begin{bmatrix}
0 & 1 \\ \nicefrac{g}{l} - \nicefrac{k}{m} & 0
\end{bmatrix},
\quad
B_2 =
\begin{bmatrix}
0 \\ \nicefrac{1}{ml^2}
\end{bmatrix},
\quad
c_2 =
\begin{bmatrix}
0 \\ \nicefrac{kd}{ml}
\end{bmatrix},\\
\mathcal D_1 &=
\{ (x,u) \ | \
x_1\!\leq\! \nicefrac{d}{l},
x_{\mathrm{min}}\!\leq x\!\leq\!x_{\mathrm{max}}, u_{\mathrm{min}}\!\leq\! u\! \leq \!u_{\mathrm{max}} \},\\
\mathcal D_2 &=
\{ (x,u) \ |
\ x_1\! >\!\nicefrac{d}{l}, x_{\mathrm{min}}\! \leq\! x\! \leq\! x_{\mathrm{max}},
 u_{\mathrm{min}}\! \leq u\! \leq\! u_{\mathrm{max}} \}
\end{align*}
The key difference between the two examples is that the cart-wall system makes and breaks contact instantaneously while the pendulum with an elastic wall allows variable contact duration, thus resulting in different mode sequences.
Due to this difference, we use a control invariant set as terminal set in this experiment 2, avoiding the terminal mode ``gets stuck in the wall''.

\begin{figure}
  \centering
  \includegraphics[width=0.75\linewidth]{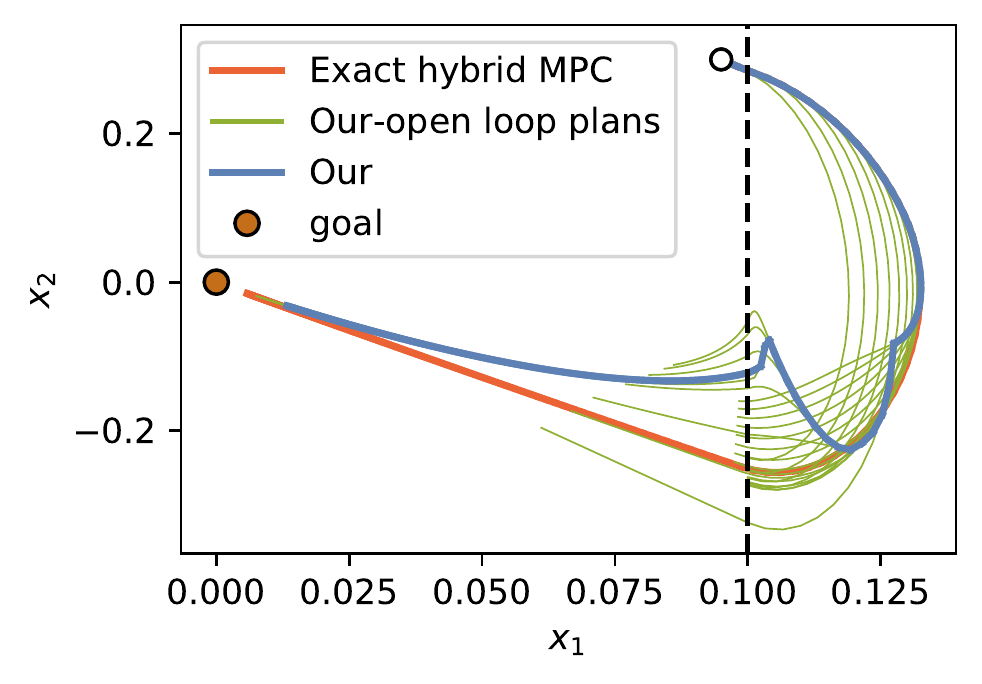}
\caption{
  Example trajectory with MPC open-loop plans for the pendulum with an elastic wall. 
}
\label{fig:e3:traj}
\end{figure}

Similar to the previous experiment, we execute Alg.\ref{alg:hmpc} to obtain the \thisalg learning-controller and collected data samples. They are plotted in Fig.~\ref{fig:e3}.
An example closed-loop trajectory in Fig.~\ref{fig:e3:traj}.
We can see that, compared with the previous example, the \thisalg trajectories differ from the optimal solution obtained by running MIP to full optimality (no time limit).
This is caused by the fact that variable-during contact results in suboptimal mode sequences being reused by Alg.~\ref{alg:hmpc}. Again, this may be either improved by Alg.~\ref{alg:imprv} or a shift-started mode sequence in HMPC.

\subsection{Computational complexity reduction}
\begin{figure}
\centering
   \includegraphics[width=.75\linewidth]{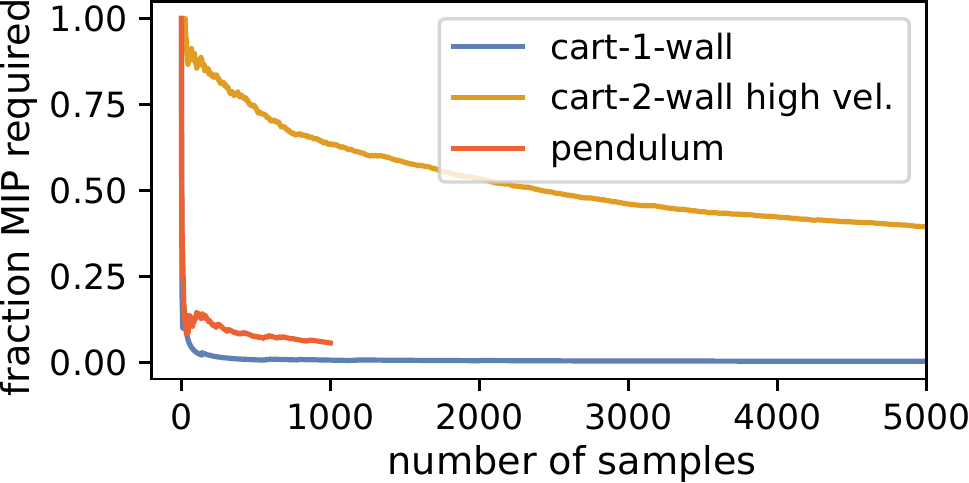}
  \vspace*{-0.3em}
  \caption{ Percentage of MIP runs during the execution of \thisalg in different environments.
    This empirically validates Proposition \ref{thm:imprv}. A preliminary wall-clock comparison is given in the appendix.
  }
\label{fig:percent}
\end{figure}
We now discuss how \thisalg in Alg.~\ref{alg:hmpc} can greatly speed up the HMPC solution method. The results here empirically validate the theory in Sec.~\ref{sec:theory}.
We compare the solution time of \thisalg in Alg.~\ref{alg:hmpc} with that of the exact MIQP solution of hybrid MPC \eqref{eq:ocp}\eqref{eq:bigm} and report the relative run-time (time of Alg.~\ref{alg:hmpc})$/$(time of exact MIQP)$\%$).
This comparison is shown in Fig.~\ref{fig:percent} for all experiments.

We observe that, as Alg.~\ref{alg:hmpc} collects more data samples, only a small fraction of mixed-integer optimization needs to be solved by branch\&bound.
As there is a big gap in computation time between branch\&bound and continuous programs, the proposed scheme drastically improve computational efficiency. The actual wall clock time may vary greatly depending on the specific learning and MPC code implementation. Nonetheless, we include a wall clock-time comparison in Appendix~\ref{sec:app:wallclock}.
\repl{This experiment shows that we have thusfar answered our question posed in Section.~1, of what we are able to \textit{learn} from HMPC runs in the same environment.}{
This experiment demonstrates the answer to our question posed in Section.~1, of what can be \textit{learned} from HMPC runs in the same environment.
}

\section{Related work}
Explicit MPC [\cite{bemporad2000piecewise, bemporad2002explicit, tondel2003algorithm}] seeks to offload online computation to offline.
Its insight is that
1)~the optimal control law of an OCP with quadratic objective linear constraints is piecewise affine state-feedback, \ie
\(u(x) = F^i x + G^i, x\in \mathcal{R}_i\), where \(\mathcal{R}_i\) are convex polyhedra (referred to as critical regions) and
2)~there may be exponentially many such regions need to be computed and stored.
In the context of hybrid systems,
different regions correspond to different (switching) mode sequences, \eg in-contact or not-in-contact.
 It is then possible to carry out the mode sequence enumeration offline and store the state-feedback-affine policy.
 Despite its elegant theory, 
 explicit MPC may incur \repl{formidable}{considerable} complexity in 1)~both online look-up and offline computation and 2)~storage of the optimal partitions.
 Many approximate MPC algorithms, \eg \cite{zeilinger2011real, domahidi2011}, trade off optimality for fast runtime and feasibility. The goal is to avoid exhaustively computing, storing, and comparing regions.


 Motivated by this challenge, many approaches have been proposed in the literature to solve such hybrid control problems.
 For example, \cite{mordatch2012discovery} suggest the use of soft contact and dynamics models to smooth the dynamics, thus allowing the use of gradient-based numerical optimization.
 The work of \cite{posa2013direct} exploits homotopy methods to compute solutions that gradually converge from relaxed to accurate ones.

Closely related to this work, the authors of \cite{marcucci2017approximate} use insights from MPC stability proof to simplify explicit MPC from critical regions to inner approximations of feasible regions, resulting in significant computational saving.
Our approach further simplifies the region computation and storage  by the non-parametric learning algorithm, sidestepping polyhedron operations completely.

Very recently, a few works show promises in applying machine learning to hybrid control.
\cite{deits2018lvis} experimented with value function and policy learning in hybrid systems. Our method in this work can also serve as an efficient oracle for policy training.
This work shares similar features with \cite{hogan2018reactive} in reducing computational cost by learning mode sequences. However, we exploit the non-forgetting property of non-parametric methods to perform online learning, as well as the geometric structure of hybrid MPC solutions to avoid exhaustive offline MIP runs.

There is a thread of works studying the topic of  ``learning to branch and bound.'' [\cite{he2014learning,khalil2016learning,rachelson2010combining}]. We share the common {thread} of using \emph{function approximation}.
\repl{
  However, those learning methods are not motivated by the structure of the control solution --- they are mostly \textit{black-box} approaches.}{
  However, those learning methods are mostly \textit{black-box} approaches.
}
In contrast, our insight is rooted in the understanding of HMPC solution structure.

\section{Discussion}
\label{sec:discuss}
This paper presents a new approach that learns from prior experiences to accelerate MIP for hybrid MPC.
Our choice of non-parametric learning algorithms is motivated by the theoretical understanding of (H)MPC solutions.
For example, neural network training is parametric and cannot easily deal with wrong points near the hybrid mode boundaries. It also looses \textit{stability guarantees} of (H)MPC.
In contrast, our method does not need training and preserves the stability of the (H)MPC.
It achieves \textit{speed-up virtually for free} by using adaptive samples while is easy-to-implement.

It is known that if absolute optimality is the goal instead of feasibility which was the aim of this paper, the boundaries of \textit{optimal mode sequence partition} may be quadratic rather than affine (cf. \cite{borrelli2017predictive}). While \thisalg still works, an interesting extension is to generalize the notion of (weighted) Euclidean distance, possibly via kernelized nearest neighbor.

One potential application of \thisalg is to serve as an efficient supervised learning oracle for training an imitation learning policy $\pi : x \mapsto u$, such as in \cite{deits2018lvis,hogan2018reactive}. Their computational cost of training may be greatly reduced by adopting \thisalg.

\begin{ack}
We would like to thank Majid Khadiv and Brahayam Ponton for their helpful discussions, Marin Vlastelica Pogancic for the help in producing Fig.~\ref{fig:ex1_traj}, Vincent Berenz and Lidia Pavel for their help in equipment and facility support.
\end{ack}

\bibliography{mioc}

\begin{thebibliography}{25}
\providecommand{\natexlab}[1]{#1}
\providecommand{\url}[1]{\texttt{#1}}
\providecommand{\urlprefix}{URL }
\expandafter\ifx\csname urlstyle\endcsname\relax
  \providecommand{\doi}[1]{doi:\discretionary{}{}{}#1}\else
  \providecommand{\doi}{doi:\discretionary{}{}{}\begingroup
  \urlstyle{rm}\Url}\fi

\bibitem[{Bemporad et~al.(2000)Bemporad, Borrelli, and
  Morari}]{bemporad2000piecewise}
Bemporad, A., Borrelli, F., and Morari, M. (2000).
\newblock Piecewise linear optimal controllers for hybrid systems.
\newblock In \emph{Proceedings of the 2000 American Control Conference. ACC},
  volume~2, 1190--1194. IEEE.

\bibitem[{Bemporad et~al.(2002)Bemporad, Morari, Dua, and
  Pistikopoulos}]{bemporad2002explicit}
Bemporad, A., Morari, M., Dua, V., and Pistikopoulos, E.N. (2002).
\newblock The explicit linear quadratic regulator for constrained systems.
\newblock \emph{Automatica}, 38(1), 3--20.

\bibitem[{Bertsimas et~al.(2016)Bertsimas, King, Mazumder
  et~al.}]{bertsimas2016best}
Bertsimas, D., King, A., Mazumder, R., et~al. (2016).
\newblock Best subset selection via a modern optimization lens.
\newblock \emph{The annals of statistics}, 44(2), 813--852.

\bibitem[{Borrelli et~al.(2017)Borrelli, Bemporad, and
  Morari}]{borrelli2017predictive}
Borrelli, F., Bemporad, A., and Morari, M. (2017).
\newblock \emph{Predictive control for linear and hybrid systems}.
\newblock Cambridge University Press.

\bibitem[{Cover et~al.(1967)Cover, Hart et~al.}]{cover1967nearest}
Cover, T.M., Hart, P., et~al. (1967).
\newblock Nearest neighbor pattern classification.
\newblock \emph{IEEE Transactions on Information Theory}, 13(1), 21--27.

\bibitem[{Deits et~al.(2018)Deits, Koolen, and Tedrake}]{deits2018lvis}
Deits, R., Koolen, T., and Tedrake, R. (2018).
\newblock {LVIS}: Learning from value function intervals for contact-aware
  robot controllers.
\newblock \emph{arXiv preprint arXiv:1809.05802}.

\bibitem[{Diehl et~al.(2002)Diehl, Bock, Schl{\"o}der, Findeisen, Nagy, and
  Allg{\"o}wer}]{diehl2002real}
Diehl, M., Bock, H.G., Schl{\"o}der, J.P., Findeisen, R., Nagy, Z., and
  Allg{\"o}wer, F. (2002).
\newblock Real-time optimization and nonlinear model predictive control of
  processes governed by differential-algebraic equations.
\newblock \emph{Journal of Process Control}, 12(4), 577--585.

\bibitem[{Domahidi et~al.(2011)Domahidi, Zeilinger, Morari, and
  Jones}]{domahidi2011}
Domahidi, A., Zeilinger, M.N., Morari, M., and Jones, C.N. (2011).
\newblock {Learning a Feasible and Stabilizing Explicit Model Predictive
  Control Law by Robust Optimization}.
\newblock \emph{Conference on Decision and Control}, 513--519.

\bibitem[{Ferreau et~al.(2008)Ferreau, Bock, and Diehl}]{ferreau2008online}
Ferreau, H.J., Bock, H.G., and Diehl, M. (2008).
\newblock An online active set strategy to overcome the limitations of explicit
  mpc.
\newblock \emph{International Journal of Robust and Nonlinear Control:
  IFAC-Affiliated Journal}, 18(8), 816--830.

\bibitem[{He et~al.(2014)He, Daume~III, and Eisner}]{he2014learning}
He, H., Daume~III, H., and Eisner, J.M. (2014).
\newblock Learning to search in branch and bound algorithms.
\newblock In \emph{Advances in neural information processing systems},
  3293--3301.

\bibitem[{Hogan et~al.(2018)Hogan, Grau, and Rodriguez}]{hogan2018reactive}
Hogan, F.R., Grau, E.R., and Rodriguez, A. (2018).
\newblock Reactive planar manipulation with convex hybrid mpc.
\newblock In \emph{2018 IEEE International Conference on Robotics and
  Automation (ICRA)}, 247--253. IEEE.

\bibitem[{Houska et~al.(2011)Houska, Ferreau, and Diehl}]{houska2011acado}
Houska, B., Ferreau, H.J., and Diehl, M. (2011).
\newblock {ACADO} toolkit--an open-source framework for automatic control and
  dynamic optimization.
\newblock \emph{Optimal Control Applications and Methods}, 32(3), 298--312.

\bibitem[{Jones et~al.(2006)Jones, Grieder, and
  Rakovi{\'c}}]{jones2006logarithmic}
Jones, C.N., Grieder, P., and Rakovi{\'c}, S.V. (2006).
\newblock A logarithmic-time solution to the point location problem for
  parametric linear programming.
\newblock \emph{Automatica}, 42(12), 2215--2218.

\bibitem[{Khalil et~al.(2016)Khalil, Le~Bodic, Song, Nemhauser, and
  Dilkina}]{khalil2016learning}
Khalil, E.B., Le~Bodic, P., Song, L., Nemhauser, G., and Dilkina, B. (2016).
\newblock Learning to branch in mixed integer programming.
\newblock In \emph{Thirtieth AAAI Conference on Artificial Intelligence}.

\bibitem[{Marcucci et~al.(2017)Marcucci, Deits, Gabiccini, Bicchi, and
  Tedrake}]{marcucci2017approximate}
Marcucci, T., Deits, R., Gabiccini, M., Bicchi, A., and Tedrake, R. (2017).
\newblock Approximate hybrid model predictive control for multi-contact push
  recovery in complex environments.
\newblock In \emph{IEEE-RAS International Conference on Humanoid Robotics},
  (Humanoids), 31--38. IEEE.

\bibitem[{Marcucci and Tedrake(2019)}]{marcucci2019mixed}
Marcucci, T. and Tedrake, R. (2019).
\newblock Mixed-integer formulations for optimal control of piecewise-affine
  systems.
\newblock In \emph{Proceedings of the Intl.\ Conference on Hybrid Systems:
  Computation and Control}, 230--239. ACM.

\bibitem[{Mordatch et~al.(2012)Mordatch, Todorov, and
  Popovi{\'c}}]{mordatch2012discovery}
Mordatch, I., Todorov, E., and Popovi{\'c}, Z. (2012).
\newblock Discovery of complex behaviors through contact-invariant
  optimization.
\newblock \emph{ACM Transactions on Graphics}, 31(4).

\bibitem[{Omohundro(1989)}]{Omohundro89fiveballtree}
Omohundro, S.M. (1989).
\newblock Five balltree construction algorithms.
\newblock Technical report, UC Berkeley.

\bibitem[{Posa and Tedrake(2013)}]{posa2013direct}
Posa, M. and Tedrake, R. (2013).
\newblock Direct trajectory optimization of rigid body dynamical systems
  through contact.
\newblock In \emph{Algorithmic Foundations of Robotics X}, 527--542. Springer.

\bibitem[{Rachelson et~al.(2010)Rachelson, Abbes, and
  Diemer}]{rachelson2010combining}
Rachelson, E., Abbes, A.B., and Diemer, S. (2010).
\newblock Combining mixed integer programming and supervised learning for fast
  re-planning.
\newblock In \emph{International Conference on Tools with Artificial
  Intelligence}, 63--70. IEEE.

\bibitem[{Rawlings and Mayne(2009)}]{rawlings2009model}
Rawlings, J.B. and Mayne, D.Q. (2009).
\newblock \emph{Model predictive control: Theory and design}.
\newblock Nob Hill Pub., Wisconsin.

\bibitem[{Richalet et~al.(1978)Richalet, Rault, Testud, and
  Papon}]{richalet1978model}
Richalet, J., Rault, A., Testud, J., and Papon, J. (1978).
\newblock Model predictive heuristic control.
\newblock \emph{Automatica (Journal of IFAC)}, 14(5), 413--428.

\bibitem[{Sontag(1981)}]{sontag1981nonlinear}
Sontag, E. (1981).
\newblock Nonlinear regulation: The piecewise linear approach.
\newblock \emph{IEEE Transactions on Automatic Control}, 26(2), 346--358.

\bibitem[{T{\o}Ndel et~al.(2003)T{\o}Ndel, Johansen, and
  Bemporad}]{tondel2003algorithm}
T{\o}Ndel, P., Johansen, T.A., and Bemporad, A. (2003).
\newblock An algorithm for multi-parametric quadratic programming and explicit
  mpc solutions.
\newblock \emph{Automatica}, 39(3), 489--497.

\bibitem[{Zeilinger et~al.(2011)Zeilinger, Jones, and
  Morari}]{zeilinger2011real}
Zeilinger, M.N., Jones, C.N., and Morari, M. (2011).
\newblock Real-time suboptimal model predictive control using a combination of
  explicit {MPC} and online optimization.
\newblock \emph{IEEE Trans.\ on Automatic Control}, 56(7), 1524--1534.

\end{thebibliography}








\appendix

\section{Wall-clock comparison}\label{sec:app:wallclock}
This section corresponds to the computation time benchmark in Fig.~\ref{fig:percent}.
We show the results using a straightforward implementation of \thisalg in Python Sklearn library.

\begin{tabular}[]{@{}lll@{}}
  Number of OCPs & LNMS runtime (s) & Regular MIP \\\hline \Tstrut
10 & \textbf{4.07} & 69.02\\
100 & \textbf{69.28} & 156.25\\
500 & \textbf{249.53} & 1186.90\\
\end{tabular}

We generate OCPs (with different initial states).
There is a big gap in solving speed between MIP and continuous programs.
Our method is significantly faster ($\sim2-17X$ speed-up).
Speed-ups are expected values across all runs.
Note, that this is done without code optimization on our side 
while MIP is solved using Gurobi -- further putting LNMS in a disadvantage.

\end{document}


\maketitle

\appendix
\section{Preliminary wall-clock comparison}
This section corresponds to the computation time benchmark in fig.5. Depending on the specific learning
and MPC implementation, the wall-clock results may vary greatly. Here are preliminary results.

\begin{tabular}[]{@{}lll@{}}
\toprule
Number of OCPs & LNMS runtime (s) & Regular MIP\tabularnewline
\midrule
\endhead
10 & \textbf{4.07} & 69.02\tabularnewline
100 & \textbf{69.28} & 156.25\tabularnewline
500 & \textbf{249.53} & 1186.90\tabularnewline
\bottomrule
\end{tabular}

We used the setting in the elastic-wall-pendulum experiment described in
Sec.4. We generate OCPs (with different initial states) and solve them
using the MPC formulation reported in the paper. There is a big gap in solving speed between MIP and continuous programs,
as we have already discussed in the paper. Our method is significantly
faster (\textasciitilde{} 2-17X speed-up). We did not cherry pick:
speed-ups are across all runs of experiment 1 (cart-wall) and another
(unreported) experiment. This is even done without code optimization and
python nearest neighbor (NN) implementation (while MIP is solved using
the fastest commercial solver. This further puts LNMS in disadvantage). In
practical MPC, this advantage can be improved drastically.

\section{Cart wall environment: equations}
Th cart with one wall environment is described by the following PWA dynamics:
\begin{equation}
\label{eq:db_int}
\begin{aligned}
\begin{cases}
x_1^+ = x_1 +  x_2\ \Delta t, \quad x_2^+ = x_2 + \dfrac{u}{m}\ \Delta t,  & \text{if } x\in \mathcal{C}_1\\
x_1^+ = x_1, \quad x_2^+ = -\epsilon x_2 , & \text{if } x\in \mathcal{C}_2
\end{cases}
\end{aligned}
\end{equation}
where $m$ is the mass constant and $\epsilon$ is the rebounce factor.
It could be thought of an actuated version of a bouncing ball---a classical hybrid system.
 $\mathcal{C}_1 = \{x_1 +  x_2\ \Delta t < x_{\text{wall}}\}$ denotes the state space where the dynamics is the double integrator and
 $\mathcal{C}_2 = \{x_1 +  x_2\ \Delta t \ge x_{\text{wall}}\}$ denotes the state space where the contact with the wall happens.
The PWA for cart with two walls are a trivial extension of the above with a third constraint set.

\section{Pendulum with elastic wall: equations}
The environment is adapted from~\cite{marcucci2017approximate}.
\begin{equation}
\dot x =
\begin{cases}
A_1 x + B_1 u, \quad \text{if} \quad (x, u) \in \mathcal D_1, \\
A_2 x + B_2 u + c_2, \quad \text{if} \quad (x, u) \in \mathcal D_2,
\end{cases}
\end{equation}

with

\begin{equation}
A_1 =
\begin{bmatrix}
0 & 1 \\ g/l & 0
\end{bmatrix},
\quad
B_1 =
\begin{bmatrix}
0 \\ 1/(ml^2)
\end{bmatrix},
\quad
A_2 =
\begin{bmatrix}
0 & 1 \\ g/l - k/m & 0
\end{bmatrix},
\quad
B_2 =
\begin{bmatrix}
0 \\ 1/(ml^2)
\end{bmatrix},
\quad
c_2 =
\begin{bmatrix}
0 \\ kd/(ml)
\end{bmatrix},
\end{equation}

\begin{equation}
\mathcal D_1 =
\{ (x,u) \ | \
x_1 \leq d/l, \
x_{\mathrm{min}} \leq x \leq x_{\mathrm{max}}, \
u_{\mathrm{min}} \leq u \leq u_{\mathrm{max}} \},
\end{equation}

\begin{equation}
\mathcal D_2 =
\{ (x,u) \ |
\ x_1 > d/l, \
x_{\mathrm{min}} \leq x \leq x_{\mathrm{max}},
\ u_{\mathrm{min}} \leq u \leq u_{\mathrm{max}} \}
\end{equation}

%
%
%
%

\begin{definition}[Critical Region]
A Critical Region is the set of parameters x for which the same set $A$ of constraints is active at the optimum. Given $x_0$,
$$
\mathcal{CR}(x_0) := \{ x\ |\ A(x) = A(x_0) \}.
$$
In the context of hybrid systems, critical regions are not polyhedral. They may have curved boundaries. Exact hybrid MPC exhaustively compare overlapped regions to determine which one defines the optimal PWA control law.
\end{definition}

\subsection{Policy training}
We learn a policy as described in Sec.~\ref{sec:mthd_policy}
To this end we used two approaches.
An naive approach that we denote as naive imitation learning (behavior cloning), where the policy is trained using randomly sampled points in state-space and labeled by Algorithm~\ref{alg:hmpc}.
move
A second approach uses DAgger which is described by Alg.~\ref{alg:dagger}
Since we learned a deterministic policy, to ensure that we eventually converge to the optimal trajectory we added Gaussian noise during the policy rollout.
Our implementation is outlined in Algorithm~\ref{alg:dagger}.
\begin{algorithm}[H]
  \caption{\thisalg: DAgger}
  \label{alg:dagger}
  \begin{algorithmic}[1]
    \STATE Initialize starting policy $\pi^*$
    \LOOP
    \STATE Rollout $N$ trajectories $\tau$ with policy $\pi^*+\mathcal{N}(0, \sigma)$
    \STATE Relabel the trajectories $\tau$ using MIQP and add to dataset $\mathcal{D}$
    \STATE Train new policy $\pi$ on dataset $\mathcal{D}$
    \STATE Set $\pi^*$ to $\pi$
    \ENDLOOP
  \end{algorithmic}
\end{algorithm}

In practice, Algorithm~\ref{alg:dagger} can be made more efficient in that the training set for the new policy is a subset of the collected dataset  $\mathcal{D}$ with a preference towards recent data to ensure continuous improvement.

We modeled the policies for both the naive imitation and DAgger approaches as neural networks with ReLU activations in the hidden layers and optimized them using the Adam optimization algorithm.
The relabeling was done using MPC with horizon 10.
The results are summarized in Figure~\ref{fig:imitation} for experiment 1.

\bibliography{mioc}